\documentclass[10pt]{amsart}

\usepackage[margin=1in]{geometry} 
\usepackage{amsmath,amsthm,amssymb, bbm, dsfont}
\usepackage{graphicx}
\usepackage{hyperref}
\usepackage{adjustbox}
\usepackage[table,xcdraw]{xcolor}
\usepackage{listings}
\usepackage{amsmath}
\usepackage{amssymb}
\usepackage{bm}
\usepackage{graphicx}
\usepackage{psfrag}
\usepackage{color}
\usepackage{hyperref}
\hypersetup{colorlinks=true, linkcolor=blue, citecolor=magenta, urlcolor=wine}
\usepackage{url}
\usepackage{algpseudocode}
\usepackage{fancyhdr}
\usepackage{mathtools}
\usepackage{tikz-cd}
\usepackage{xy}
\input xy
\xyoption{all}
\usepackage{stmaryrd}
\usepackage{calrsfs}

\newtheorem{thm}{Theorem}[section]
\newtheorem{prop}[thm]{Proposition}

\newtheorem{lemma}[thm]{Lemma}
\newtheorem{corll}[thm]{Corollary}
\theoremstyle{definition}
\newtheorem{defn}[thm]{Definition}
\newtheorem{rem}[thm]{Remark}
\newtheorem{example}[thm]{Example}

\numberwithin{equation}{section}

\newcommand{\tor}{\textsf{T}}

\newcommand{\cc}{\mathbb{C}}

\newcommand{\nn}{\mathbb{N}}
\newcommand{\pp}{\mathbb{P}}
\newcommand{\ppp}{\mathcal{P}}
\newcommand{\qq}{\mathbb{Q}}

\newcommand{\zz}{\mathbb{Z}}
\newcommand{\oo}{\mathcal{O}}

\newcommand{\spec}{\text{Spec}}
\newcommand{\vp}{\textsf{v}}
\newcommand{\rank}{\text{rank }}

\newcommand{\lcm}{\text{lcm}}

\providecommand\ldb{\llbracket}
\providecommand\rdb{\rrbracket}

\newcommand\gp{\mathsf{gp}}

\newcommand{\lc}{\text{lc}}

\newcommand{\supp}{\text{supp }}
\newcommand{\pow}{\text{pow}}

\newcommand{\ord}{\text{ord }}

\newcommand{\at}{\mathcal{A}}


\begin{document}
 
% --------------------------------------------------------------
%                         Start here
% --------------------------------------------------------------
\mbox{}
\title{On the atomicity of one-dimensional monoid algebras}

\author{Ishan Panpaliya}
\address{Department of Mathematics\\Seattle University\\Seattle, WA 98122, US}
\email{ipanpaliya@seattleu.edu}

\date{\today}

\begin{abstract}
	The ascending chain condition on principal ideals (ACCP) is almost always complementary to atomicity within integral domains: in fact, Cohn initially stated that these two conditions were equivalent. This assertion has been shown to be false, however most counterexamples require technical algebraic constructions. In 2017, Gotti conjectured that for every $q$ in the set $S := ((0, 1) \cap \mathbb{Q}) \setminus {\mathbb{N}}^{-1}_{> 1}$, atomicity ascends from the exponentially cyclic Puiseux monoid $M_q$ to its monoid algebra over the field of rationals. If this conjecture were true, it would provide an extremely wide class of atomic domains of Krull dimension one not satisfying the ACCP, and so would be perhaps the simplest possible such examples. Bu et al. recently proved that the monoid algebra $\mathbb{Q} \left[M_{3/4} \right]$ is atomic, marking the first progress towards settling this conjecture. We strengthen this result and prove that $\mathbb{Q}[M_q]$ is atomic for all $q \in S$ having an odd denominator.
\end{abstract}

\bigskip
%%%%%%%%%%%%%%%%
\maketitle

\section{Introduction}
\label{sec: intro}
A unital ring is said to be \emph{atomic} if every nonunit element of its multiplicative monoid has a factorization into atoms (i.e. irreducible elements). Atomicity has been systematically studied in the setting of commutative rings and monoids since it was introduced by Cohn in 1968~\cite{pC68}. Specifically, atomicity has been investigated in Pr\"ufer domains~ \cite{rH16, CGH23, rH21}, rings of integer-valued polynomials~\cite{ACCS94}, Noetherian and Cohen-Kaplansky rings~\cite{AAV01, CA10}, and semigroup rings~\cite{fG20, GL23}. There is also an extensive literature characterizing, in particular, the atomicity of Puiseux monoids as well as their corresponding monoid algebras, and they have been shown to admit surprisingly rich atomic structure (see~\cite{CGG21, CJMM24, GL23, GR24} and references therein). Tringali has also recently given a general characterization by evaluating orderings on monoids~\cite{sT23}. A comprehensive survey on the history and open problems regarding atomicity can be found in~\cite{CG24}. 

An integral domain satisfies the \emph{ascending chain condition on principal ideals} (ACCP) if every ascending sequence of its principal ideals stabilizes. In~\cite{pC68}, Cohn infamously stated without proof that every integral domain that is atomic also satisfies the ACCP. This assertion has been shown to be false, although both properties are closely related to one another. In general, the ACCP is a strictly stronger condition than atomicity, but both are equivalent if further conditions are imposed. For example, it is known that atomicity and the ACCP are equivalent within the class of irreducible divisor finite (IDF) monoids. In addition, for a torsion-free monoid $M$, the condition that every submonoid of $M$ is atomic is equivalent to that of $M$ satisfying the ACCP~\cite{GV23}. 

As opposed to the Noetherian condition, the ACCP deals only with principal ideals. Hence, it admits an equivalent characterization which can naturally be extended to the more general setting of commutative and cancellative monoids: for every $x$ in the multiplicative monoid of an integral domain $R$,
given a sequence $X := (x_n)_{n \in \nn}$ of elements in $R^*$ such that 
    \[
    \cdots \mid_{R} \ x_4 \ \mid_{R} \ x_3 \ \mid_{R} \ x_2 \ \mid_R \ x_1 := x,
    \]
there exists $t \in {\nn}$ such that $x_t \sim_R x_s$ for all $s \in \nn_{\geq t}$. Such a characterization motivates the following generalization of the ACCP given in~\cite{GL24}. A monoid $M$ is said to satisfy the \emph{almost ACCP} if for every finite, nonempty subset $S$ of $M$, there exists a common divisor $d \in M$ of $S$ such that for some $s \in S$, the element $\frac{s}{d}$ satisfies the ACCP in $M$. Both the ACCP and the almost ACCP are intimately linked to atomicity. In fact, the following diagram holds, with counterexamples known for the reverse of each implication~\cite{GL24}: $\text{ACCP} \implies \text{Almost ACCP} \implies \text{Atomic}$. The ACCP is a simple necessary condition for an integral domain to be Noetherian and its failure indicates that a domain lacks the unique factorization property. 

The first atomic domain failing to satisfy the ACCP was constructed by Grams in 1974~\cite{aG74}. It was constructed by localizing the one-dimensional monoid algebra $\qq[G]$ at its set of polynomials with nonzero order, where the Grams monoid, $G$, is defined as follows (here $p_n$ denotes the $n^{\text{\tiny th}}$ odd prime):
\[
G := \left\langle 
\frac{1}{2^n p_n}: n \in \nn
\right\rangle.
\]
It is still an open problem as to whether $\qq[G]$ itself is an atomic domain.
The following polynomial ring was conjectured to be another example by Cohn and only proven so by Zaks in 1982~\cite{aZ82}, where $F$ is an arbitrary field and $U_i := \frac{X_0 Y}{X_i Z^i}$ for all $i \in \nn_0$:
\[
F\left[
\{X_i\}^{\infty}_{i=0}, Y, Z, \{U_i\}^{\infty}_{i=0}
\right].
\]
In 1998, Roitman provided a quite technical example using Rees Algebras, which came as a by-product of his work settling the ascent of atomicity to polynomial extensions~\cite{mR93}. More recently, examples have been constructed by Boynton and Coykendall~\cite{BC19}, and by Gotti and Li~\cite{GL23, GL24}, who were able to provide two examples of monoid algebras whose monoids have rank at least two.

Monoid algebras of Puiseux monoids over the field of rationals are an obvious class of candidates for simple atomic integral domains failing to satisfy the ACCP: they are natural extensions of polynomial rings whose monoids have rank equal to $\nn_0$, and, moreover, are integral domains of dimension one. Such monoid algebras have been studied by Gotti, who studied the boundedness, finite factorization, and ACCP~\cite{fG22}, and by Benelmekki and El Baghdadi, who focused on the IDF property~\cite{BB22}. It is well known that the ACCP ascends from a monoid to its monoid algebra, regardless of the ring of coefficients. Therefore, finding such an example amounts to proving the atomicity of the monoid algebra $\qq[M]$, for some atomic Puiseux monoid $M$ which does not satisfy the ACCP. Towards this end, in~\cite{fG22} Gotti conjectured that atomicity ascends from the class of exponentially cyclic monoids to their correponding monoid algebras over $\qq$. 

Exponentially (or multiplicatively cyclic) Puiseux monoids were introduced by Gotti and Gotti in 2016~\cite{GG18}. They constitute a broad class of rank one monoids which are parametrized by the positive rationals. For any $q \in \qq_{> 0}$, the exponentially cyclic monoid $M_{q}$ is defined as follows:
\[
M_{q} := \left \langle q^n : n \in \nn_0 \right \rangle.
\]
The atomicity, ACCP, finite factorization, and boundedness, as well as related arithmetic invariants of exponentially cyclic monoids have been investigated by Albizu-Campos, Bringas, and Polo~\cite{ABP21}; Chapman, Gotti, and Gotti~\cite{CGG20}; and Polo~\cite{hP20}. With regards to atomic structure, it is known that $\at(M_q) = \{q^n : q \in \nn_0 \}$ when $q \notin \nn^{\pm 1}$. When $q \in \nn$, the set $\at(M_q)$ consists of only one element, and when $q \in \nn_{> 1}^{-1}$, the monoid $M_q$ is antimatter. As $\min \at(M_q) = 1$ when $q > 1$, in this case the monoid $M_q$ satisfies both the finite factorization and bounded factorization properties (and so the ACCP as well). On the other hand, when $q \in ((0,1) \cap \qq) \setminus \nn^{-1}_{> 1}$, the monoid $M_q$ is atomic yet simultaneously does not satisfy the ACCP. In a recent paper, Gotti and Li have proved that exponentially cyclic monoids which are not antimatter also satisfy the almost ACCP~\cite[Example~3.8]{GL23}.

Gotti's conjecture is relatively uninteresting when $q \in \qq_{\geq 1} \cup \nn_{> 1}^{- 1}$, as the monoid $M_q$ and the monoid algebra $\qq[M_q]$ both satisfy the ACCP or both are not atomic. However, proving the atomicity of $\qq[M_q]$ for $q \in ((0, 1) \cap \qq) \setminus \nn_{> 1}^{-1}$ would provide an exceptionally broad family of one-dimensional monoid algebras which are atomic but do not satisfy the ACCP. Although the conjecture was posed in 2017, this problem has remained completely open; that is, until recently when Bu et al. proved that $\qq[M_{3/4}]$ is atomic~\cite{BGLZ24}. This provides the first piece of evidence towards resolving the conjecture hypothesized in~\cite{fG22}. Here, we strengthen this result and prove that $\qq[M_q]$ is atomic for all $q \in ((0, 1) \cap \qq) \setminus \nn^{-1}_{> 1}$ having an odd denominator.

An outline of the paper is as follows. In Section~\ref{sec:background}, we first introduce the notation we use, as well as all necessary mathematical background needed for the paper. Then, we recall results regarding the factorization of cyclotomic polynomials and composed polynomials of the form $f(x^r)$, where $f(x)$ is an irreducible polynomial in $\zz[x]$ and $r$ is an exponent in $\nn$. Finally, we introduce families of polynomials which we require in later sections of the paper. In Section~\ref{sec: SS and CDS}, we introduce a generalization of splitting sequences, which were introduced in~\cite{BGLZ24}, as well as several definitions associated with them. The most technical part of our paper is Section~\ref{sec: lambdas}, in which we prove bounds on the lengths of splitting sequences of certain families of polynomials introduced in Section~\ref{sec:background}. This section also develops methods from algebraic number theory, which we believe will be of use in tackling the more general questions regarding the ascent of atomicity, finite factorization, and bounded factorization to monoid algebras. In Section~\ref{sec: atomic exp cyc}, we prove our main result, which is that the monoid algebra $\qq[M_q]$ is atomic for all rationals $q$ in the set $((0,1) \cap \qq) \setminus \nn_{> 1}^{-1}$ having an odd denominator. This provides the first infinite class of monoid algebras of Krull dimension one that are atomic but do not satisfy the ACCP, and makes significant progress towards resolving Gotti's conjecture.

%%%%%%%%%%%%%%%%
\bigskip

\section{Background}
\label{sec:background}

%%%%%%%%%%%%%%%%
\subsection{General Notation}
\smallskip

Here we adopt the convention that the natural numbers refers to the set of positive integers. The notation $\ldb a, b \rdb$ will be used to denote the discrete interval between two integers $a, b \in \zz$; that is, $\ldb a, b \rdb := \{x \in \zz: a \leq x \leq b\}$. Given a subset $S$ of the real numbers, and a real number $r$, the notation $S_{\geq r}$ refers to the set of all elements of $S$ greater than or equal to $r$. Analogously, we may also write $S_{> r}$, $S_{\leq r}$, or $S_{< r}$. Let $q \in \qq_{> 0}$. We denote by $\textsf{n}(q)$ and \textsf{d}($q$) the unique pair of natural numbers such that $\gcd(\textsf{n}(q), \textsf{d}(q)) = 1$ and $q = \frac{\textsf{n}(q)}{\textsf{d}(q)}$. We let $\pp$ denote the set of prime numbers. For any $p \in \pp$, we let $\vp_p: \qq^* \rightarrow \zz$ denote the $p$-adic valuation map. Given a natural number $n \in \nn_{\geq 2}$, the function $\omega(n)$ counts the number of distinct prime divisors of $n$ and $\Omega(n)$ the number of prime divisors of $n$ with multiplicity (we adopt here the convention that $\omega(1) = \Omega(1) = 0$). Furthermore, we let $\pi(n)$ be the set of all distinct prime divisors of $n$ and $\Pi(n)$ the collection of all prime divisors of $n$ with multiplicity so that $|\pi(n)| = \omega(n)$ and $|\Pi(n)| = \Omega(n)$. Now let $P$ be a finite set of primes. We set $P_{\nn} := \{ n \in \nn_{\geq 2}: \pi(n) \subseteq P \}$ so that any natural number in $P_{\nn}$ has prime divisors only in $P$. Given any complex number $c \in \cc$ and any pair of natural numbers $(n,r) \in \nn \times \nn_{\geq 2}$, the notation $c^{r^{- n}}$ denotes an arbitrary ${(r^n)}^{\text{\tiny th}}$ root of $c$. Suppose that we have a finite sequence $A$ and a (possibly infinite) sequence $B$. We denote by $A \circ B$ the concatenation of the sequences $A$ and $B$.

\medskip
%%%%%%%%%%%%%%%%%%%

\subsection{Commutative Monoids, Factorizations, and Ideals}

A \emph{monoid} is a semigroup with an identity element. We assume that every monoid under consideration here is both commutative and cancellative. Let $M$ be a monoid written additively. The set of $\emph{invertible}$ elements of $M$, also referred to as the $\emph{units}$ of $M$, is denoted by $M^{\times}$. The monoid $M$ is called \emph{torsion-free} provided that, for elements $a,b \in M$, the equality $na = nb$ implies that $a = b$, for any index $n \in \nn$. We denote by $\tor(M)$ the subgroup of $M^{\times}$ consisting of its torsion elements. Given elements $a, b \in M$, we say that $a$ is \emph{associate} to $b$ if there exists $u \in M^{\times}$ such that $a + u = b$. In this case, we write $a \sim_M b$. The \emph{rank} of $M$, denoted by $\text{rank} \, M$, is the dimension of the vector space $\qq \otimes_\zz \gp(M)$ over~$\qq$, where $\gp(M)$ is the \emph{Grothendieck group} of $M$. A special class of rank one monoids is that of additive submonoids of $\qq_{\geq 0}$, which are called $\emph{Puiseux monoids}$. If $M$ is not a (nontrivial) group, in addition to being torsion-free and having rank one, then $M$ is isomorphic to a Puiseux monoid (see~\cite[Section~24]{lF70} and \cite[Theorem~3.12]{GGT21}). For a subset $S$ of $M$, we let $\langle S \rangle$ denote the submonoid of $M$ generated by~$S$.
\smallskip

An element $a \in M \setminus M^{\times}$ is called an $\emph{atom}$ if the fact that the equality $a = b + c$ holds, for elements $b,c \in M$, implies that either $b \in M^{\times}$ or $c \in M^{\times}$. We denote the set of all atoms of $M$ by $\at(M)$. An element of $M$ is called \emph{atomic} provided that it can be expressed as a sum of atoms and/or units. Therefore, the condition that $M$ is atomic is equivalent to that of each element of $M$ being atomic. 

\medskip
%%%%%%%%%%%%%%%%

\subsection{Ideals, Number Fields, and Monoid Algebras} 

A commutative ring with identity is called an \emph{integral domain} if the product of any two nonzero elements is nonzero. Let $R$ be an integral domain. We let $R^*$ denote the monoid consisting of all nonzero elements of $R$ under the multiplicative operation of $R$. $R^*$ is called the \emph{multiplicative monoid} of $R$. If for any $a,b \in R^*$, the equivalence relation $a \sim_{R^*} b$ holds, we write $a \sim_{R} b$ instead. Similarly, if the divisibility relation $a \mid_{R^*} b$ holds, we write $a \mid_{R} b$ instead.
\smallskip

Let $R$ be an integral domain. A subgroup, $I$, of $R$ is called an \emph{ideal} if for any $(x, r) \in I \times R$, the inclusion $x r \in I$ holds. If $I$ is of the form $b R$, for some element $b \in R$, then $I$ is called a \emph{principal ideal} and $I$ is called a \emph{prime ideal} if the equality $a b \in I$ implies that either $a \in I$ or $b \in I$, for any elements $a,b \in R$. We denote by $\spec(R)$ the set of all prime ideals of $R$. A sequence of ideals $\mathcal{I} := (I_n)_{n \in \nn}$ is said to be \emph{ascending} if for all $n \in \nn$, the inclusion $I_n \subseteq I_{n+1}$ holds. $\mathcal{I}$ is said to $\emph{stabilize}$ if there exists $N \in \nn$ such that for all $n \in N_{\geq N}$, the equality $I_n = I_N$ holds.
\smallskip

Let $M$ be a Puiseux monoid. As is standard, the notation $R[x;M]$ denotes the commutative ring consisting of all polynomial expressions in the indeterminate $x$ with coefficients in $R$ and exponents in $M$. As $R[x; M]$ is an $R$-module, we call $R[x; M]$ the \emph{monoid algebra} of $M$ 
over $R$. Following Gilmer~\cite{rG84}, we shorten the notation to $R[M]$ from here on out. Since $M$ is torsion-free, the monoid algebra $R[M]$ is an integral domain~\cite[Theorem~8.1]{rG84}. Moreover, since $M$ is a subset of $\qq$, it is endowed with a natural ordering $(\preceq)$. Then, it is not hard to see that any nonzero polynomial $f(x) \in R[M]$ can be written in the following manner, for an index $k \in \nn$, nonzero coefficients $r_1, \ldots, r_k \in R^*$, and distinct exponents $m_1, \ldots, m_k \in M$ such that $m_1 \prec \dots \prec m_k$:
\begin{equation} \label{eq:poly in mon alg}
	f(x) = r_1 x^{m_1} + \cdots + r_k x^{m_k}.
\end{equation}
With $f(x)$ written as above, we define the \emph{support} of $f(x)$ to be the set $\supp f := \{ m_1, \ldots, m_k \}$. The polynomial $f(x)$ is said to be \emph{ACCP-supported} if at least one element of $\supp f$ satisfies the ACCP in $M$. The exponents $\ord f := m_1$ and $\deg f := m_k$ are referred to as the \emph{order} and the \emph{degree} of $f(x)$, respectively. In the case that $R$ is a GCD-domain, we define the \emph{content} of $f(x)$ to be the quantity $\textsf{c}(f) := \gcd_R(r_1,\ldots,r_k)$.
\smallskip

We now prove that atomicity in $\qq[M_q]$ is implied by atomicity in $\zz[M_q]$; that is, in a certain sense we can translate factorizations between the two monoid algebras. This ability is crucial in the proof of Proposition~\ref{prop: (general) num irr. div is limited}, as it allows us to invoke the fact that for any polynomial $f(x) \in \zz[x]$, there are only finitely many possibilities for the leading coefficient of any of its divisors.

\begin{prop}\label{prop: Z[M] atomic = Q[M] atomic}
    Let $q$ be a positive rational and let $M_q$ be the exponentially cyclic Puiseux monoid parametrized by $q$. If the monoid algebra $\zz[M_q]$ is atomic, then the monoid algebra $\qq[M_q]$ is atomic.
\end{prop}
\begin{proof}
    We prove the lemma directly. Let $f(x)$ be an arbitrary nonunit polynomial in the monoid algebra $\qq[M_q]$. It is enough to show that $f(x)$ has an atomic factorization in $\qq[M_q]$. As the coefficients of $f(x)$ are rationals, we can take $\gamma \in \zz$ such that $\gamma f(x) \in \zz[M_q]$. By assumption, $\gamma f(x)$ has an atomic factorization in $\zz[M_q]$ and so we may write $\gamma f(x) = \mu \prod_{i = 1}^m{b_i(x)}$, for an index $m \in \nn$, irreducible nonconstant polynomials $b_1(x),...,b_m(x) \in \zz[M_q]$, and a constant polynomial $\mu \in \zz$. Then the factorization $f(x) = \frac{\mu}{\gamma} \prod_{i = 1}^m{b_i(x)}$ holds in $\qq[M_q]$. We show that this is an atomic factorization. Since $\frac{\mu}{\gamma} \in \qq$, this simply amounts to showing that for each $i \in \ldb 1, m \rdb$, the inclusion $b_i(x) \in \at(\qq[M_q])$ holds: this is a direct consequence of~\cite[Theorem~3.3]{fG20}.
\end{proof} 

\begin{rem}\label{rem: atoms of z[1/p^k N] = z[1/p^l N]}
    Let $r$ be any nonzero positive rational and let $R$ be an integral domain. As $\nn_0$ is isomorphic to any nonzero multiple of itself, there is an $R$-algebra isomorphism $\varphi: R[x] \rightarrow R[r \nn_0]$ defined by $\varphi(f(x)) = f(x^r)$ for an arbitrary polynomial $f(x) \in R[x]$. Therefore, factorizations may be translated between these isomorphic monoid algebras after applying $\varphi$ or $\varphi^{-1}$.
\end{rem}

The \emph{Krull dimension} of $R$ is the supremum of the lengths of all descending chains of its prime ideals (by inclusion). We denote the Krull dimension of $R$ by $\dim R$. If $R$ is Noetherian, then it follows from~\cite[Corollary~2]{jO88} that $\dim R[M] = \dim R + \text{rank} \, M$. Therefore, a monoid algebra that is an integral domain is one-dimensional if and only if it is the algebra of a rank one torsion-free monoid over a field. More information regarding monoid algebras and progress on the ascent of algebraic properties from the pair $(R,M)$ to the algebra $R[M]$, up to 1984, can be found in Gilmer's book~\cite{rG84}.

\medskip
%%%%%%%%%%%%%%%%

\subsection{Number Fields and Cyclotomic Polynomials} 

Let $\alpha$ be a complex number. We say that $\alpha$ is \emph{integral} over $\zz$ if $\alpha$ is the root of a monic polynomial whose coefficients lie in $\zz$.
Let $K$ be a field containing $\qq$. If $K$ has a finite basis as a vector space over $\qq$, then $K$ is called a \emph{number field}. Suppose now that $F$ is any arbitrary number field. We denote by $\mathcal{O}_F$ the \emph{ring of integers} of $F$; that is, the set of all elements of $F$ integral over $\zz$.
\smallskip

Let $n$ be a natural number. We denote by $\zeta_n$ a primitive $n^{\text{\tiny th}}$ root of unity. We use the notation $\Phi_n(x)$ to denote the unique monic polynomial in $\qq[x]$ whose root set is precisely the set of primitive $n^{\text{\tiny th}}$ roots of unity (i.e. the $n^{\text{\tiny th}}$ cyclotomic polynomial). It is well known that given a prime $p$ and an index $y \in \nn$, the equality $\Phi_{y}(x^p) = \Phi_{yp}(x)$ holds when $p$ divides $y$. In the case that $p$ does not divide $y$, the polynomial $\Phi_y(x^p)$ splits in $\zz[x]$ as $\Phi_y(x^p) = \Phi_y(x) \Phi_{yp}(x)$. From this, we can deduce the following elementary identity.

\begin{lemma}\label{lemma: cyclotomicFactors}
    Let $y$ be an index in $\nn$ and let $r \in \nn_{\geq 2}$ be a squarefree positive integer. Then the following atomic factorization holds in $\zz[x]$:
    \[
    \Phi_y(x^r) = 
    \prod_{d \mid_{\zz} \frac{r}{\gcd(y,r)}}{
    \Phi_{y \cdot d \cdot \gcd(y,r)}(x).
    }
    \]
\end{lemma}

We also recall another useful cyclotomic identity.
\begin{lemma}\label{lemma: factor-x**n-1}
    Let $n \in \nn$. Then the polynomial $x^n - 1$ factors in $\zz[x]$ as $x^n - 1 = \prod_{d \mid_{\zz} n}{\Phi_d(x)}$.
\end{lemma}

%%%%%%%%%%%%%%%
\medskip

\subsection{Composed Polynomials}

It is a result of Guersenzvaig~\cite[Theorems~4.3 and~4.4]{nG09} that for any odd prime $p$ and any irreducible polynomial $f(x) \in \zz[x]$ of degree $d \in \nn$, if the composed polynomial $f(x^p)$ is reducible over $\zz$, then there exist polynomials $S_0(x), S_1(x), \ldots, S_{p-1}(x) \in \zz[x]$ such that
\[
    f(x^p) \sim_{\zz[x]} \text{det}
    \left(
    \begin{bmatrix}
        S_0(x^p) & x S_1(x^p) & \cdots & x^{p-1} S_{p-1}(x^p)\\
        x^{p-1}S_{p-1}(x^p) & S_0(x^p) & \cdots & x^{p-2} S_{p-2}(x^p)\\
        \vdots & \vdots & \ddots & \vdots\\
        x S_1(x^p) & x^2 s_2(x^p) & \cdots & S_0(x^p)
    \end{bmatrix}
    \right)
\] 
and the polynomial $P(x) := \sum_{i=0}^{p-1}{S_i(x^p) x^i}$ is an irreducible polynomial that divides $f(x^p)$ in $\zz[x]$ with leading coefficient associate to either ${\pm \lc(f)}^{1/p}$ in $\zz$. With a little bit of extra work, we can deduce the following lemma.

\begin{lemma} \label{lemma: |lc| <= |lc|}
    Let $f(x)$ be an irreducible polynomial in $\zz[x]$ and let $r \in \nn_{\geq 3}$ be a squarefree odd positive integer. Letting $p_{\text{max}}$ be the maximal prime dividing $r$, if $f(x^r)$ is reducible, then the magnitude of the leading coefficient of any of its irreducible divisors in $\zz[x]$ is at most ${|\lc(f)|}^{\frac{p_{\text{max}} - 1}{p_{\text{max}}}}$.
\end{lemma}
Let $P'(x) := \frac{f(x^p)}{P(x)}$ so that $P'(x)$ is a divisor of $f(x^p)$ in $\zz[x]$. It was also proven in~\cite{nG09} that every irreducible divisor of $P'(x)$ in $\zz[x]$ has degree
\[
    [\qq(\alpha) : \qq(\alpha) \cap \qq(\zeta_p)] \cdot (p-1),
\]
where $\alpha$ is any root of $f(x)$. Since $[\qq(\alpha) \cap \qq(\zeta_p) : \qq] \leq p - 1$, and the minimal polynomial of $\alpha$ over $\qq$ has degree $d$, the following (in)equalities hold:
\[
    [\qq(\alpha) : \qq(\alpha) \cap \qq(\zeta_p)] = \frac{[\qq(\alpha) : \qq]}{[\qq(\alpha) \cap \qq(\zeta_p) : \qq]} \geq  \frac{d}{p-1},
\]
and so any irreducible divisor of $P'(x)$ has degree at least $d$. In particular, $f(x^p)$ has at most $p$ irreducible divisors in $\zz[x]$ and so we obtain the following lemma.

\begin{lemma} \label{lemma: <= r irr. divs}
    Let $f(x) \in \zz[x]$ be an irreducible polynomial and let $r \in \nn_{\geq 3}$ be an odd squarefree positive integer. If the polynomial $f(x^r)$ is reducible, then $f(x^r)$ has at most $r$ irreducible divisors in $\zz[x]$.
\end{lemma}

We now recall a well-known result regarding the irreducibility of composed polynomials over arbitrary fields.

\begin{lemma}[Capelli's Lemma]\label{fullPolyComposeIrr}
    Let $K$ be a field and let $f(x) \in K[x]$ be an irreducible polynomial with root $\alpha \in \overline{K}$. Furthermore, let $r \in \nn_{\geq 3}$ be an odd positive integer. Then $f(x^r)$ is irreducible over $K$ if and only if $\alpha \notin {[\qq(\alpha)]}^p$ for every prime $p \mid_{\zz} r$.
\end{lemma}

We now introduce an important set of polynomials which we will consider in Sections~\ref{sec: SS and CDS} and~\ref{sec: lambdas} of the paper.

\begin{defn} \label{defn: composed div sets}
    Let $f(x)$ be an irreducible polynomial in $\zz[x]$. Fix a set $T \subsetneq \nn_0$. We define the \emph{composed divisor set} of $f(x)$ with respect to $T$ as follows:
    \[
    \mathcal{D}(f(x), T) := \{ p(x) \in \zz[x]: \exists r \in T \text{ such that } p(x) \mid_{\zz[x]} f(x^r) \}.
    \]
    Given some finite $\omega \in \nn_0$, we also define the \emph{exponent-restricted composed divisor set} of $f(x)$ with respect to $T$ as
    \[
    \mathcal{D}(f(x), T, \omega) :=
    \{ p(x) \in \zz[x]: \exists r \in T \text{ such that } \omega(r) \leq \omega \text{ and } p(x) \mid_{\zz[x]} f(x^r) \}.
    \]
    It is not difficult to see that the cardinality of the set $\mathcal{D}(f(x), T, \omega)$ is finite.
\end{defn}

%%%%%%%%%%%%%%%%
\bigskip

\section{On Splitting Sequences}
\label{sec: SS and CDS}

Splitting sequences were introduced by Bu et al. in~\cite{BGLZ24}, where they were used to investigate the atomicity and ACCP in the monoid algebra $\qq[M_{3/4}]$. In this section, our aim is to expand upon (and slightly tweak) the definition of these splitting sequences so as to admit general exponents and permit the investigation of factorizations in arbitrary rank one monoid algebras. Our definition allows us to analyze chains of composed polynomials of the form $f(x^r)$, where the exponent $r$ can be an arbitrary natural number.

\begin{defn} \label{defn: splitSeqDef}
    Let $f(x)$ be an irreducible polynomial in $\zz[x]$. We will construct a sequence of irreducible polynomials $F := (f_n(x))_{n \in \nn_0}$ whose first term is $f(x)$. To do so, we will need a sequence of exponents: let $E_F := (e_n)_{n \in \nn}$ be a sequence of positive integers. We refer to $E_F$ as the \emph{exponent sequence} of $F$.
    To construct $F$, first set $f_0(x) := f(x)$. Now suppose that we have already defined the polynomials $f_0(x),\ldots,f_{k-1}(x) \in F$, for some $k \in \nn$. Then we simply set $f_{k}(x)$ to be any irreducible divisor of $f_{k-1}(x^{e_k})$ in $\zz[x]$. We say that $F$ is an $E_F$ \emph{splitting sequence} of $f(x)$. 
    Note that an arbitrary polynomial $f(x) \in \zz[x]$ may have multiple distinct splitting sequences. In addition, we may sometimes only specify a finite number of terms in the exponent sequence $E_F$, in which case it is tacitly assumed that a finite number of terms in $F$ are well-defined while their cofinite complement is irrelevant in the specific context.
\end{defn}

\begin{defn} \label{def: bin str def}
    We now introduce (and slightly modify) the definition of a \emph{binary string} introduced in~\cite{BGLZ24}. 
    Let $f(x)$ be as in Definition~\ref{defn: splitSeqDef} and
    let $F := (f_n(x))_{n \in \nn_0}$ be an $ E_F$ splitting sequence of $f(x)$, where $E_F := (e_n)_{n \in \nn}$. We then associate to each polynomial in $F$, besides $f(x)$, a digit from the alphabet \{``L", ``S"\}, corresponding to the words ``Lift" and ``Split". The sequence of these digits forms an infinite sequence, which we refer to as the \emph{binary string} of $F$. We denote this sequence by $\mathcal{B}(F)$ and index its $n^{\text{\tiny th}}$ term by $\mathcal{B}_n(F)$. For the construction, if $f_n(x) = f_{n-1}(x^{e_{n}})$, then set $\mathcal{B}_{n}(F) := ``L"$; otherwise $f_{n-1}(x^{e_n})$ is reducible in $\zz[x]$ and so we set $\mathcal{B}_n(F) := \text{``S"}$. 
\end{defn}

\begin{defn}
    Let $f(x)$ be as in Definition~\ref{defn: splitSeqDef}. We denote by $\mathcal{S}( f(x))$ the \emph{splitting set} of $f(x)$; that is, the set of all splitting sequences of $f(x)$. If we fix an index $m \in \nn$ and a set $T := \{ t_1, \ldots, t_m \}$ of distinct natural numbers, then we denote by $\mathcal{S}(f(x), T)$ the subset of $\mathcal{S}(f(x))$ whose members have exponent sequences with exponents in $T$. For any $F \in \mathcal{S}(f(x))$, we define $\textsf{N}_{\text{spl}}(F)$ as the number of splits at the beginning of $F$. Equivalently, $\textsf{N}_{\text{spl}}(F)$ is defined to be the minimal index in $\nn$ such that the equality $f_{\textsf{N}_{\text{spl}}(F)+1}(x) = f_{\textsf{N}_{\text{spl}}(F)}(x^{e_{\textsf{N}_{\text{spl}}(F)+1}})$ holds.
    We then define
    \[
    \Lambda(f(x), T) := \sup_{F \in \mathcal{S}(f(x), T)}{\textsf{N}_{\text{spl}}(F)}
    \]
    so that $\Lambda(f(x), T)$ returns the supremum of the number of initial splits over all splitting sequences of $f(x)$ with exponent sequence whose elements are contained in $T$.
\end{defn}

To obtain a clearer understanding of splitting sequences, we present a few different examples.

\begin{example} \hfill
    \begin{enumerate}
        \item [(1)] Let $f(x) := x$. Let $F := (f_n(x))_{n \in \nn_0}$ be any splitting sequence of $f(x)$ with exponent sequence $E_F := (e_n)_{n \in \nn}$, where $e_n > 1$ for all $n \in \nn$. Then $f_n(x) \in \{ \pm x \}$ for all $n \in \nn$ and
        $\mathcal{B}(F) = (\text{``S"})_{n \in \nn}$.
        \smallskip
        
        \item [(2)] Let $f(x) := \Phi_1(x)$ and, setting $p_n$ to be the $n^{\text{\tiny th}}$ prime for all $n \in \nn$, let $E_F := (p_n)_{n \in \nn}$ be an exponent sequence. Then $F := \Phi_1(x) \circ (\Phi_{\prod_{i=1}^n{p_n}}(x))_{n \in \nn_0}$ is an $E_F$ splitting sequence of $f(x)$ with $\mathcal{B}(F) = (\text{``S"})_{n \in \nn}$.
        \smallskip

        \item[(3)] Let $f(x) := 1 + 2x + 3x^2 + x^3$ and let $E_F := (3)_{n \in \nn}$ be an exponent sequence. Then, letting $g(x) := 1 - x + x^3$, the sequence $F := f(x) \circ (g(x^{3^{n-1}}))_{n \in \nn}$ is an $E_F$ splitting sequence of $f(x)$ with $\mathcal{B}(F) = \text{``S"} \circ (\text{``L"})_{n \in \nn}$.
    \end{enumerate}
\end{example}

\begin{lemma}\label{lemma: irr div appears in some splitting seq}
    Let $f(x)$ be an irreducible nonconstant polynomial in $\zz[x]$. For any irreducible polynomial $b(x) \in \zz[x]$ not associate to $f(x)$ such that there exists $r_b \in \nn_{\geq 2}$ for which $b(x) \mid_{\zz[x]} f(x^{r_b})$, for any $k \in \nn$ and any sequence $\sigma := (\sigma_n)_{n \in \ldb 1, k \rdb}$ of natural numbers greater than one satisfying that $\prod_{i=1}^k{\sigma_i} = r_b$, there exists an $E_F := \sigma$ splitting sequence, $F := (f_n(x))_{n \in \nn_0}$, of $f(x)$ such that $f_k(x) = b(x)$. 
\end{lemma}
\begin{proof}
    We prove the claim by induction on $k$. For the base case, we take $k = 1$. Suppose now that $d(x)$ is an irreducible polynomial in $\zz[x]$ such that there exists $r_d \in \nn_{\geq 2}$ such that $d(x) \mid_{\zz[x]} f(x^{r_d})$. The equality $\prod_{i=1}^{1}{\sigma_i} = r_d$ implies that $\sigma = (r_d)$, and so proving the base case amounts to proving that there exists an $(r_d)$ splitting sequence $F := (f_n(x))_{n \in \nn_0}$ of $f(x)$ such that $f_1(x) = d(x)$. This is trivial. For the inductive step, suppose that there exists $K \in \nn$ such that the lemma holds for all $k \in \ldb 1, K \rdb$. Let $d(x)$ be as in the base case, and suppose now that $\sigma := (\sigma_n)_{n \in \ldb 1, K+1 \rdb}$ is a sequence of natural numbers greater than one satisfying that $\prod_{i=1}^{K+1}{\sigma_i} = r_d$. First note that the divisibility relation $d(x) \mid_{\zz[x]} f({(x^{\prod_{i=1}^K{\sigma_i}})}^{\sigma_{K+1}})$ holds in $\zz[x]$. Let $s(x)$ be the unique irreducible divisor of $f(x^{\prod_{i=1}^K{\sigma_i}})$ in $\zz[x]$ such that $d(x) \mid_{\zz[x]} s(x^{\sigma_{K+1}})$. The inductive hypothesis guarantees the existence of a splitting sequence $F := (f_n(x))_{n \in \nn_0}$ with $E_F := (\sigma_n)_{n \in \ldb 1, K \rdb}$ such that $f_K(x) = s(x)$. Then the fact that $d(x)$ is irreducible in $\zz[x]$ and that $d(x) \mid_{\zz[x]} s(x^{\sigma_{K+1}})$ ensures that we may set the ${(K+1)}^{\text{\tiny th}}$ exponent in $E_F$ to be $\sigma_{K+1}$ and $f_{K+1}(x)$ to be $d(x)$, proving the lemma for $k = K+1$.
\end{proof}

\bigskip
%%%%%%%%%%%%%%

\section{Bounding Lengths of Splitting Sequences in $\zz[x]$} \label{sec: lambdas}
In order to prove our main theorem in Section~\ref{sec: atomic exp cyc}, we will need to show that ACCP-supported polynomials in $\zz[M_q]$ satisfy the ACCP, a condition stronger than atomicity itself. In order to do so, we must necessarily gain control over the factorization of polynomials in $\zz[x]$. In particular, given a fixed set of primes $P$ and an irreducible polynomial $f(x) \in \zz[x]$, we must provide absolute bounds on $\Lambda(b(x), P)$ across all polynomials $b(x)$ in the composed divisor set $\mathcal{D}(f(x), P)$ (aside from a few exceptional cases, that is).

We first recall a useful result which may be deduced from Capelli's Lemma.
\begin{lemma}\label{lemma: IrrTowerK}
    Let $K$ be a field. Given an irreducible polynomial $f(x) \in K[x]$ and an odd positive integer $r \in \nn_{\geq 3}$, if $f(x^r)$ is irreducible over $K$, then $f(x^{r^k})$ is also irreducible over $K$ for any $k \in \nn_{\geq 2}$.
\end{lemma}

Lemma~\ref{lemma: IrrTowerK} will be pivotal in the rest of the proof as it provides insight into the structure of splitting sequences. In particular, we obtain the following easy corollary.

\begin{corll}\label{corll: SSS...}
    Let $f(x) \in \zz[x]$ be an irreducible polynomial and let $r \in \nn_{\geq 3}$ be an odd positive integer. If $F$ is an $(r)_{n \in \nn}$ splitting sequence of $f(x)$ that has $T$ splits, for some $T \in \nn$, then $\mathcal{B}(F)$ starts with the sequence $(\text{``S"})_{n \in \ldb 1, T \rdb}$.
\end{corll}

We now introduce the function $\Lambda^*$, which we will use in order to bound the number of initial splits in any splitting sequence starting at $f(x)$ (that is, the quantity $\Lambda(f(x), T)$), where $f(x)$ is some fixed irreducible polynomial in $\zz[x]$. The main goal of this section is to prove the inequality $\Lambda(f(x),T) \leq \Lambda^*(f(x), T)$ for as wide a class of polynomials in $\zz[x]$ as possible, focusing primarily on the case when $T$ is a finite set of primes.

\begin{defn} \label{defn: lambda*}
%%FIX K_{f, i, P}
    Let $f(x)$ be an irreducible polynomial in $\zz[x]$ with $\lc(f) \in \{\pm 1\}$ and degree $d_f \in \nn$. Let $m$ be an index in $\nn$ and $P := \{p_1, \ldots, p_m \}$ a set of primes. For each $i \in \ldb 1, d_f \rdb$, let $\alpha_i$ be a distinct root of $f(x)$ so that $A_f := \left\{ \alpha_{f,1},\ldots,\alpha_{f,d_f} \right\}$ is the root set of $f(x)$; let $K_{f,i}$ denote the number field $\qq(\alpha_{f,i})( \zeta_{p_1}, \ldots, \zeta_{p_m})$; let $T_{f,i} := \rank \oo^{\times}_{K_{f, i}}$; and let $G_{f, i} := \left\{g(f,i,1), \ldots, g(f,i,T_{f,i}) \right\}$ denote a minimal set of generators of  $\oo^{\times}_{K_{f, i}} / \tor\left(\oo^{\times}_{K_{f, i}}\right)$ that is fixed from here on out. Then we define
    \[
    \Lambda^{\star}(f(x), P) := \max_{i \in \ldb 1, d_f \rdb}
    \begin{dcases}
        \sum_{\substack{\mathfrak{p} \in \spec\left(\oo_{K_{f,i}}\right)\\ \mathfrak{p}| \alpha_{f,i} \oo_{K_{f,i}}}}{
            \vp_{\mathfrak{p}}(\alpha_{f,i} \oo_{K_{f,i}})
        },
        & \alpha_{f,i} \in \oo^*_{K_{f,i}} \setminus \oo^{\times}_{K_{f,i}},\\
        \sum_{p \in P}{
            \min_{
                \substack{k \in \ldb 1, T_{f,i} \rdb \\ \vp_{g(f,i,k)}(\alpha_{f,i}) \neq  0}
            }{
                \vp_{p}(\vp_{g(f,i,k)}(\alpha_{f,i}))
            }
        },
        &\alpha_{f,i} \in \oo^{\times}_{K_{f,i}} \setminus \tor\left(\oo^{\times}_{K_{f,i}}\right),\\
        \infty, &
        \text{otherwise}.
    \end{dcases}
    \]
\end{defn}

\begin{lemma} \label{lemma: Finite SS Lambda(f,p) lc = -1,1}
    Let $m$ be an index in $\nn$ and $P := \{ p_1, \ldots, p_m \}$ a set of primes. Let $f(x) \in \zz[x] \setminus (\{\pm x\} \cup \{ \pm \Phi_y(x): y \in \nn \})$ be an irreducible polynomial with $\lc(f) \in \{\pm 1\}$. Assume also that the terminology of Definition~\ref{defn: lambda*} associated to $f(x)$ and $P$ holds. Then the quantity $\Lambda(f(x), P)$ is at most $\Lambda^{\star}(f(x), P)$.
\end{lemma}
\begin{proof}
    To prove the lemma, we assume by way of contradiction that $f(x)$ and $P$ are as in the statement of the lemma, but that $\Lambda(f(x), P)$ is not bounded above by $\Lambda^{\star}(f(x), P)$. We first replace $f(x)$ by its monic associate in $\zz[x]$ (note that this substitution is harmless as the sign of $f(x)$ has no effect on any divisibility relation in $\zz[x]$). To shorten notation, we write instead $\Lambda^{\star}$ instead of $\Lambda^{\star}(f(x),P)$ hereafter (we may do so as $f(x)$ and $P$ are fixed).
    By the initial assumptions of the proof, we can find a splitting sequence $F := (f_n(x))_{n \in \nn_0}$ of $f(x)$ with exponent sequence $E := (e_n)_{n \in \nn}$ and $\textsf{N}_{\text{spl}}(F) \geq \Lambda^{\star} + 1$ such that for all $n \in \nn$, the inclusion $e_n \in P$ holds. Since $f(x)$ is monic, for each $n \in \nn_0$, the leading coefficient of $f_n(x)$ must be in the set $\{\pm 1\}$. As we did before, we replace each polynomial in $F$ with its monic associate. Now fix a root $\beta$ of $f_{\Lambda^{\star}+1}(x)$. The fact that the divisibility relation $f_{n}(x) \mid_{\zz[x]} f_{n-1}(x^{e_{n}})$ holds for every $n \in \nn$ implies that $\beta^{\prod_{n=1}^{\Lambda^{\star}+1}{e_n}}$ is a root of $f(x)$. Then there exists $i_{\beta} \in \ldb 1, d_f \rdb$ such that $\alpha_{f, i_{\beta}} = \beta^{\prod_{n=1}^{\Lambda^{\star}+1}{e_n}}$. 
    \smallskip

    \noindent \textsc{Claim.} $\beta$ is contained in the number field $K_{f, i_{\beta}}$.
    \smallskip

    \noindent \textsc{Proof of Claim.} We prove the following statement by induction: for every $n \in \ldb 0, \Lambda^{\star} + 1 \rdb$, the element $\beta^{\prod_{i=n+1}^{\Lambda^{\star}+1}{e_n}}$ is contained in $K_{f, i_{\beta}}$. It is clear that this proves the claim. For the base case, we take $n = 0$. The proof is trivial by the definition of $K_{f, i_{\beta}}$. Suppose now that there exists $n' \in \ldb 0, \Lambda^{\star} \rdb$ such that the inclusion $\beta^{\prod_{i=n+1}^{\Lambda^{\star}+1}{e_n}} \in K_{f, i_{\beta}}$ holds for all $n \in \ldb 0, n' \rdb$. Since $n'$ is less than $\Lambda^{\star} + 1$, the fact that $\textsf{N}_{\text{spl}}(F) \geq \Lambda^{\star}+1$ guarantees that the polynomial $f_{n'}(x^{e_{n'+1}})$ is reducible over $\zz$. Specifically, $f_{n'+1}(x)$ is an irreducible proper divisor of $f_{n'}(x^{e_{n'+1}})$ with root $\beta^{\prod_{i=n'+2}^{\Lambda^{\star} + 1}{e_n}}$. Then Capelli's Lemma guarantees that some ${(e_{n'+1})}^{\text{\tiny th}}$ root of $\beta^{\prod_{i=n'+1}^{\Lambda^{\star}+1}{e_n}}$ exists in the field $\qq\left(\beta^{\prod_{i=n'+1}^{\Lambda^{\star}+1}{e_n}}\right)$,
    whence we can conclude immediately that there exists $\ell \in \ldb 0, e_{n'+1} - 1 \rdb$ such that the inclusion 
    \[
    \zeta^{\ell}_{e_{n'+1}} \cdot \beta^{\prod_{i=n'+2}^{\Lambda^{\star} + 1}{e_n}} \in \qq\left(\beta^{\prod_{i=n'+1}^{\Lambda^{\star}+1}{e_n}}\right)
    \]
    holds. However, by the inductive hypothesis, $\beta^{\prod_{i=n'+1}^{\Lambda^{\star}+1}{e_n}}$ is contained within $K_{f, i_{\beta}}$ and so it is a direct consequence that $\zeta^{\ell}_{e_{n'+1}} \cdot \beta^{\prod_{i=n'+2}^{\Lambda^{\star} + 1}{e_n}}$ is contained in $K_{f, i_{\beta}}$. Then the proof of the claim follows from the fact that $e_{n'+1}$ is a prime in $P$ and for all $p \in P$, the number field $K_{f, i_{\beta}}$ contains a primitive $p^{\text{\tiny th}}$ root of unity.
    \smallskip
    
    Since all polynomials in $F$ are monic, it is clear that both $\beta$ and $\alpha_{f, i_{\beta}}$ are contained within the integral closure of $\zz$. This observation, in tandem with the previous claim, allows us to conclude that both $\beta$ and $\alpha_{f, i_{\beta}}$ are elements of $\oo_{K_{f, i_{\beta}}}$. We now case on whether $\alpha_{i_{\beta}}$ is zero, a nonzero nonunit, a nontorsion unit, or a torsion element in the ring $\oo_{K_{f, i_{\beta}}}$. However, since $f(x)$ is not a monomial, it must be the case that $\alpha_{f, i_{\beta}} \neq 0$. Furthermore, the fact that $f(x)$ is not associate to a cyclotomic polynomial ensures that $\alpha_{f, i_{\beta}} \notin \tor\left(\oo^{\times}_{K_{f, i_{\beta}}}\right)$. Therefore, we only need to contend with the following two cases.
    \smallskip

    \textsc{Case 1:} $\alpha_{f, i_{\beta}} \in {\oo^*_{K_{f, i_{\beta}}}} \setminus  \oo^{\times}_{K_{f, i_{\beta}}}$. The fact that $\alpha_{i_{\beta}}$ is a nonzero nonunit in $\oo_{K_{f, i_{\beta}}}$ guarantees that the same is true of $\beta$, whence we can conclude that there exists an index $\ell \in \nn$, prime ideals $\mathfrak{p}_1,\ldots,\mathfrak{p}_{\ell} \in \spec\left(\oo_{K_{f, i_{\beta}}}\right)$, and exponents $t_1, \ldots, t_{\ell} \in \nn$ such that the principal ideal $\beta \oo_{K_{f, i_{\beta}}}$ uniquely decomposes as $\beta \oo_{K_{f, i_{\beta}}} = \prod_{j = 1}^{\ell}{\mathfrak{p}_j^{t_j}}$. The equality $\alpha_{f, i_{\beta}} = \beta^{\prod_{n=1}^{\Lambda^{\star}+1}{e_n}}$ then ensures that $\alpha_{f, i_{\beta}} \oo_{K_{f, i_{\beta}}} = {(\beta \oo_{K_{f, i_{\beta}}})}^{\prod_{n=1}^{\Lambda^{\star}+1}{e_n}}$. The previous statements together result in the following factorization:
    \begin{equation} \label{prin. idl. fac.}
        \alpha_{f, i_{\beta}} \oo_{K_{f, i_{\beta}}} = \prod_{j = 1}^{\ell}{\mathfrak{p}_j^{t_j \cdot \prod_{n=1}^{\Lambda^{\star}+1}{e_n}}}.
    \end{equation}
    Therefore, we obtain that
    \[
    \sum_{
        \substack{
        \mathfrak{p} \in \spec\left(\oo_{K_{f,i_{\beta}}}\right)\\
        \mathfrak{p} | \alpha_{f,i_{\beta}} \oo_{K_{f,i_{\beta}}}
        }
    }{
        \vp_{\mathfrak{p}}{
        \left(
        \alpha_{f,i_{\beta}} \oo_{K_{f,i_{\beta}}}
        \right)
        }
    } = \sum_{j = 1}^{\ell}{t_j \prod_{n=1}^{\Lambda^{\star}+1}{e_n} }.
    \]
    Let $p_{\text{min}}$ denote the minimal prime in $P$. By the definition of $\Lambda^{\star}$, in tandem with the minimality of $p_{\text{min}}$, we recover the following inequality:
    \[
    \sum_{
        \substack{
            \mathfrak{p} \in \spec\left(\oo_{K_{f,i_{\beta}}}\right) \\
            \mathfrak{p} | \alpha_{f,i_{\beta}} \oo_{K_{f,i_{\beta}}}
        }
    }{
        \vp_{\mathfrak{p}}{
        \left(
        \alpha_{f,i_{\beta}} \oo_{K_{f,i_{\beta}}}
        \right)
        }
    } >
    \pow\left(
    p_{\text{min}}, \sum_{
        \substack{
            \mathfrak{p} \in \spec\left(\oo_{K_{f,i_{\beta}}}\right)\\
            \mathfrak{p} | \alpha_{f,i_{\beta}} \oo_{K_{f,i_{\beta}}}
        }
    }{
        \vp_{\mathfrak{p}}{
        \left(
        \alpha_{f,i_{\beta}} \oo_{K_{f,i_{\beta}}}
        \right)
        }
    }
    \right).
    \]
    However, this is clearly a contradiction. As such, we are forced to conclude that this case is invalid.
    \smallskip

    \textsc{Case 2:} $\alpha_{f, i_{\beta}} \in \oo^{\times}_{K_{f,  i_{\beta}}} \setminus \tor\left(\oo^{\times}_{K_{f,  i_{\beta}}}\right)$. Since for every $n \in \nn$, the exponent $e_n$ is a prime in the set $P$, the equality $\alpha_{f, i_{\beta}} = \beta^{\prod_{n=1}^{\Lambda^{\star}+1}{e_n}}$ may be rewritten as $\alpha_{f, i_{\beta}} = \beta^{p_1^{a_1} \cdots p_m^{a_m}}$ for exponents $a_1,\ldots,a_m \in \nn_0$ such that $\sum_{j=1}^{m}{a_j} = \Lambda^{\star}+1$. However, $\beta$ may also be expressed (up to a torsion element) as a product of the generators of $\oo^{\times}_{K_{f, i_{\beta}}}$; that is, the elements in the set $G_{f,i_{\beta}}$. Therefore, we may write 
    \[
    \beta  = \zeta \prod_{j=1}^{T_{f, i_{\beta}}}{ 
    {g(f, i_{\beta}, j)}^{t_j}
    },
    \]
    for some torsion element $\zeta \in \oo^{\times}_{K_{f, i_{\beta}}}$ and unique exponents $t_1,\ldots,t_{T_{f, i_{\beta}}} \in \zz$. We immediately recover the equality 
    \[
    \alpha_{f, i_{\beta}} = \zeta' \prod_{j = 1}^{T_{f, i_{\beta}}}{
    {g(f,i_{\beta}, j)}^{t_j \cdot p_1^{a_1} \cdots p_m^{a_m}}
    },
    \]
    for some multiple $\zeta'$ of $\zeta$. The fact that $\sum_{j=1}^{m}{a_j} = \Lambda^{\star}+1$ implies that for any $k \in \ldb 1, T_{f, i_{\beta}} \rdb$, the following equality holds:
    \begin{equation} \label{eq: sum v_{g_k}}
        \sum_{p \in P}{
        \vp_p(\vp_{g(f, i_{\beta}, k)}(\alpha_{f, i_{\beta}})) 
    } = 
    \begin{dcases}
        (\Lambda^{\star} + 1) + \sum_{p \in P}{
        \vp_p(\vp_{g(f, i_{\beta}, k)}(\beta))} &
        \vp_{g(f, i_{\beta}, k)}(\beta) \neq 0\\
        0, & \vp_{g(f, i_{\beta}, k)}(\beta) = 0.
    \end{dcases}
    \end{equation}
    Then the fact that $\beta$ is not a torsion element of $\oo_{K_{f, i_{\beta}}}$ guarantees that for at least one $k \in \ldb 1, T_{f, i_{\beta}} \rdb$, the $\vp_{g(f, i_{\beta}, k)}$-adic valuation of $\beta$ is nonzero, and hence at least equal to $\Lambda^{\star} + 1$ by~\eqref{eq: sum v_{g_k}}. Then taking a minimum over all $k \in \ldb 1, T_{f, i_{\beta}} \rdb$ yields that
    
    \[
    \sum_{p \in P}{
        \min_{
            \substack{k \in \ldb 1,     T_{f,i_{\beta}} \rdb\\ \vp_{g(f, i_{\beta}, k)}(\alpha_{f, i_{\beta}}) \neq  0}
        }
        {
            \vp_p(\vp_{g(f, i_{\beta}, k)}(\alpha_{f, i_{\beta}}))
        }
    }
    \geq \Lambda^{\star} + 1.
    \]
    However, this is a clear contradiction to the definition of $\Lambda^{\star}(f(x), P)$. 
    
    As both cases are invalid, the initial assumption that $\Lambda(f(x), P)$ was not bounded above by $\Lambda^{\star}(f(x), P)$ must have been false, proving the lemma.
\end{proof}

\begin{lemma}\label{lemma: Lambda(b(x),P) <= Lambda^*(f(x),P)}
    Let $m$ be an index in $\nn$ and $P := \{ p_1, \ldots, p_m \}$ a set of odd primes. Let $f(x)$ be an irreducible polynomial $f(x) \in \zz[x] \setminus (\{ \pm x \} \cup \{ \pm \Phi_y(x): y \in \nn \})$ with $\lc(f) \in \{\pm 1\}$. Assume also that the terminology of Definition~\ref{defn: lambda*} associated to $f(x)$ and $P$ holds. Then for all polynomials $b(x) \in \mathcal{D}(f(x), P)$, the quantity $\Lambda(b(x), P)$ is at most $\Lambda^{\star}(f(x), P)$.
\end{lemma}
\begin{proof}
    For each polynomial $b(x) \in \mathcal{D}(f(x), P)$, we denote by $r_{b}$ the minimal exponent in $P_{\nn}$ such that the divisibility relation $b(x) \mid_{\zz[x]} f(x^{r_b})$ holds. Now the fact that $\lc(f) \in \{ \pm 1 \}$ ensures that the same is true of all polynomials in $\mathcal{D}(f(x), P)$. Moreover, since $f(x)$ is not a monomial nor associate to a cyclotomic polynomial, the set $\mathcal{D}(f(x), P) \cap (\{ \pm x \} \cup \{ \pm \Phi_y(x): y \in \nn \})$ is empty. Then by virtue of Lemma~\ref{lemma: Finite SS Lambda(f,p) lc = -1,1}, it is sufficient to prove the statement that for all polynomials $b(x) \in \mathcal{D}(f(x), P)$, the bound $\Lambda^{\star}(b(x), P) \leq \Lambda^{\star}(f(x), P)$ holds. We prove the statement by induction on $\Omega(r_b)$; that is, the number of prime factors of $r_b$ with multiplicity. We take the base case to be $\Omega(r_b) = 0$: this is simply the statement of Lemma~\ref{lemma: Finite SS Lambda(f,p) lc = -1,1}. So now for the inductive step, assume that there exists $K \in \nn_0$ such that for all polynomials $b(x) \in \mathcal{D}(f(x), P)$ with $\Omega(r_b) \leq K$, the inequality $\Lambda^{\star}(b(x), P) \leq \Lambda^{\star}(f(x), P)$ holds. Now fix any polynomial $w(x) \in \mathcal{D}(f(x), P)$ such that $\Omega(r_w) = K+1$. Let $p \in P$ be a prime dividing $r_w$ and write $r_w = p \cdot r_w'$ for some $r_w ' \in \nn$. Then, relabel the indices of $P$ so that $p_m = p$. The divisibility relation $w(x) \mid_{\zz[x]} f(x^{p \cdot r_w '})$ implies that the roots of $w(x)$ are complex $p^{\text{\tiny th}}$ roots of the roots of $f(x^{r_w'})$. Therefore, there exists a polynomial $s(x) \in \zz[x]$ which is the unique irreducible divisor of $f(x^{r_w'})$ such that $w(x) \mid_{\zz[x]} s(x^p)$. 
    
    We will prove the bound $\Lambda^{\star}(w(x), P) \leq \Lambda^{\star}(f(x), P)$ as follows. For any index $i \in \ldb 1, d_w \rdb$, we demonstrate that the following inequality holds if $\alpha_{w,i} \in \oo^*_{K_{w,i}} \setminus \oo^{\times}_{K_{w,i}}$, where $y_i$ is an index in $\ldb 1, d_s \rdb$ such that ${(\alpha_{w, i})}^p = \alpha_{s, y_i}$:
    \begin{equation} \label{eq: first comp-wise ineql}
        \sum_{
        \substack{
            \mathfrak{q} \in \spec\left(\oo_{K_{w,i}}\right)\\
            \mathfrak{q} \mid \alpha_{w,i} \oo_{K_{w,i}}
        }
        }{
            \vp_{\mathfrak{q}}\left(\alpha_{w,i} \oo_{K_{w,i}}\right)
        } \leq 
        \sum_{
        \substack{
            \mathfrak{p} \in \spec\left(\oo_{K_{s,y_i}}\right)\\
            \mathfrak{p} \mid \alpha_{s,y_i} \oo_{K_{s,y_i}}
        }
        }{
            \vp_{\mathfrak{p}}\left(\alpha_{s,y_i} \oo_{K_{s,y_i}}\right)
        }.
    \end{equation}
    We may do so as the inductive hypothesis guarantees that the inequality $\Lambda^{\star}(s(x), P) \leq \Lambda^{\star}(f(x), P)$ holds. On the other hand, if $\alpha_{w,i} \in \oo^{\times}_{K_{w,i}} \setminus \tor\left(\oo^{\times}_{K_{w,i}}\right)$, we simply show that the inequality
    \begin{equation} \label{eq: sum vp(v-g-k) ineq}
        \sum_{p \in P}{
            \min_{\substack{k \in \ldb 1, T_{w,i} \rdb  \\ 
            \vp_{g(w, i,k)}(\alpha_{w,i}) \neq 0}}{
                \vp_{p}(\vp_{g(w, i,k)}      (\alpha_{w,i}))
            }
        } \leq \Lambda^{\star}(f(x),P)
    \end{equation}
    holds. Since $w(x) \notin \{\pm x\}  \cup  \{ \pm \Phi_y(x): y \in \nn \}$, we do not have to contend with the case when $\alpha_{w, i} = 0$ or when $\alpha_{w, i} \in \tor\left(\oo^{\times}_{K_{w,i}}\right)$. 
    
    Accordingly, fix $i' \in \ldb 1, d_w \rdb$. Since the field extension $K_{w, i'} / K_{s, y_{i'}}$ is obtained by adjoining a $p^{\text{\tiny th}}$ root of $\alpha_{s, y_{i'}}$, it is clear that the degree of the extension is either equal to one or to $p$. If $[K_{w, i'} : K_{s, y_{i'}}] = 1$, then either~\eqref{eq: first comp-wise ineql} or~\eqref{eq: sum vp(v-g-k) ineq} can be deduced easily. Hence, we can now assume that the field extension $K_{w, i'} / K_{s, y_{i'}}$ is of degree $p$. This guarantees that $\alpha_{s, y_{i'}}$ has no $p^{\text{\tiny th}}$ roots in the number field $K_{s, y_{i'}}$. We now case on whether $\alpha_{w, i'}$ is a nonzero nonunit or a nontorsion unit within the ring of integers $\oo_{K_{w, i'}}$.
    \smallskip

    \textsc{Case 1:} $\alpha_{w,i} \in \oo^*_{K_{w,i}} \setminus \oo^{\times}_{K_{w,i}}$. It is obvious that $\alpha_{s,y_{i'}}$ is also a nonzero nonunit within the ring $\oo_{K_{s, y_{i'}}}$. Let $\mathfrak{p} \in \spec\left(\oo_{K_{s, y_{i'}}}\right)$ be any prime ideal dividing $\alpha_{s, y_{i'}} \oo_{K_{s,y_{i'}}}$. The number of prime ideals, $\mathfrak{q} \in \spec \left( \oo_{K_{w,i'}} \right)$, lying over $\mathfrak{p}$ is at most $p$ since the degree of the field extension $K_{w, i'} / K_{s, y_{i'}}$ is $p$. This guarantees that the following inequality holds:
    \[
    \sum_{
    \substack{
        \mathfrak{q} \in \spec\left(\oo_{K_{w,i'}}\right)\\
        \mathfrak{q} \mid \alpha_{s, y_{i'}} \oo_{K_{w, i'}}
    }
    }{
    \vp_{\mathfrak{q}}\left(\alpha_{s, y_{i'}} \oo_{K_{w, i'}}\right)
    }
    \leq p 
    \sum_{
        \substack{\mathfrak{p} \in \spec\left(\oo_{K_{s, y_{i'}}}\right) 
        \\ \mathfrak{p} \mid \alpha_{s, y_{i'}} \oo_{K_{s, y_{i'}}}}
    }{
    \vp_{\mathfrak{p}}\left(\alpha_{s, y_{i'}}\oo_{K_{s, y_{i'}}}\right)
    }.
    \]
    The proof of the case is then a direct consequence of the factorization ${(\alpha_{w,i'} \oo_{K_{w,i'}})}^p$ = $\alpha_{s,y_{i'}}\oo_{K_{w,i'}}$.
    \smallskip

    \textsc{Case 2:} $\alpha_{w,i'} \in \oo^{\times}_{K_{w,i'}} \setminus \tor\left(\oo^{\times}_{K_{w,i'}}\right)$. For each $k \in \ldb 1, T_{w,i'} \rdb$, let $v_k := \vp_{g(w, i', k)}(\alpha_{w,i'})$ so that $\alpha_{w, i'} = \zeta \prod_{j=1}^{T_{w,i'}}{{g(w, i', j)}^{v_j}}$, for some torsion element $\zeta \in \oo_{K_{w,i'}}$. Now suppose by way of contradiction that~\eqref{eq: sum vp(v-g-k) ineq}
    does not hold. This guarantees that we can find exponents $a_1,\ldots,a_m \in \nn_0$ such that $a_1 + \cdots + a_m \geq \Lambda^{\star}(f(x),P) + 1$ and for any $k \in \ldb 1, T_{w, i'} \rdb$, we have that $v_k = v_k' \cdot p_1^{a_1} \cdots p_m^{a_m}$ for some $v_k' \in \nn_0$. It is then clear that we may express $\alpha_{w,i'}$ as follows, where $\zeta'$ is some multiple of $\zeta$ and where we choose each of the exponents $a_1,\ldots,a_m$ maximally (that is, not all of the exponents $v_1', \ldots, v_{T_{w,i'}}'$ are divisible by a prime in $P$):
    \[
    \alpha_{w,i'} = 
    {\left(
    \zeta' \prod_{j=1}^{T_{w,i'}}{
    {g(w, i', j)}^{v_j'}
    }
    \right)}^{p_1^{a_1} \cdots p_m^{a_m}}.
    \]
     Now set
    \[
    \widetilde{\alpha} := \zeta' \prod_{j=1}^{T_{w,i'}}{
    {g(w, i', j)}^{v_j'}
    }
    \]
    so that ${\widetilde{\alpha}}^{p_1^{a_1} \cdots p_m^{a_m}} = \alpha_{w,i'}$. Furthermore, set $L := 2 + a_1 + \cdots + a_m$. We now construct a sequence of exponents, $S := (s_{\ell})_{\ell \in \ldb 1, L \rdb}$, as follows. To do so, for each $i \in \ldb 1, m \rdb$, first construct the subsequence $P_i := (p \cdot p_{i+1} \cdots p_{m} \cdot  p_i^j)_{j \in \ldb 1, a_i \rdb}$. Then we define
    \[
    S := (1, p) \circ P_m \circ P_{m-1} \circ \cdots \circ P_1.
    \]
    It is not difficult to see that for every $\ell \in \ldb 1, L \rdb$, the element ${\widetilde{\alpha}}^{s_{\ell}}$ is contained in $\oo^{\times}_{K_{w, i'}}$.
    \smallskip
    
    \noindent \textsc{Claim.} $\alpha_{s, y_{i'}} \in \oo^{\times}_{K_{s,y_{i'}}}$.
    \smallskip

    \noindent \textsc{Proof of Claim.}
    The fact that $\alpha_{w,i'}$ is a unit in the ring of integers $\oo_{K_{w, i'}}$ ensures that $\alpha_{s, y_{i'}}$ also is. Therefore, $\frac{1}{\alpha_{s, y_{i'}}}$ is integral over $\zz$, whence we we can conclude that $\frac{1}{\alpha_{s, y_{i'}}}$ is contained within $\oo_{K_{s, y_{i'}}}$.
    \smallskip
    
    Observe that if there exists $\ell \in \ldb 1, L-1 \rdb$ such that $\widetilde{\alpha}^{s_{\ell}}$ is contained in $\oo^{\times}_{K_{s,y_{i'}}}$, then this implies that $\widetilde{\alpha}^{s_{\ell + z}}$ is also contained in $\oo^{\times}_{K_{s,y_{i'}}}$ for all $z \in \ldb 1, L - \ell \rdb$. This observation, in tandem with the fact that $\widetilde{\alpha}^{s_{L}} = \alpha_{s, y_{i'}} \in \oo^{\times}_{K_{s,y_{i'}}}$, guarantees that there exists a minimal index $\ell_0 \in \ldb 1, L\rdb$ such that for all $\ell \in \ldb \ell_0, L \rdb$, the element $\widetilde{\alpha}^{s_{\ell}}$ is contained within $\oo^{\times}_{K_{s,y_{i'}}}$. However by the inductive hypothesis, the inequality 
    \[
    \sum_{p \in P}{
            \min_{\substack{k \in \left\ldb 1, T_{s,y_{i'}} \right\rdb \\ \vp_{g(s,y_{i'},k)}\left(\alpha_{s,y_{i'}}\right) \neq 0}}{
                \vp_{p}\left(
                \vp_{g(s,y_{i'},k)}    \left(\alpha_{s,y_{i'}}\right)
                \right)
            }
        } \leq \Lambda^{\star}(f(x),P)
    \]
    holds. This observation, in conjunction with the fact that the sum of the exponents $a_1, \ldots, a_m$ is at least equal to $\Lambda^{\star}(f(x), P) + 1$, ensures that $\ell_0$ is at least equal to three. Set $p_1 := \frac{\ell_0 - 1}{\ell_0 - 2}$ and $p_2 := \frac{\ell_0}{\ell_0-1}$.
    \smallskip

    \textsc{Subcase 2.1:} 
    $p_1 = p_2$. It is clear that ${\widetilde{\alpha}}^{\ell_0}$ is an element of $ \oo^{\times}_{K_{s,y_{i'}}}$. At the same time, both the exclusion ${({\widetilde{\alpha}}^{\ell_0})}^{1/p_
    2} \notin \oo^{\times}_{K_{s,y_{i}}}$ and the inclusion ${({\widetilde{\alpha}}^{\ell_0})}^{1/p^2_2} \in \oo^{\times}_{K_{w,i'}}$ hold simultaneously. However, this is impossible, as this implies that the intermediate field $K_{s, y_{i'}}\left({({\widetilde{\alpha}}^{\ell_0})}^{1/p^2_2}\right)$ of the field extension $K_{w, i'} / K_{s, y_{i'}}$ has degree $p_2^2$ over $K_{s, y_{i'}}$.
    \smallskip

    \textsc{Subcase 2.2:} $p_1 \neq p_2$. It is clear that $\widetilde{\alpha}^{\ell_0}$ is an element of $\oo^{\times}_{K_{s,y_{i'}}}$. We then have that the exclusion ${\left(\widetilde{\alpha}^{\ell_0}\right)}^{1/p_2} \notin \oo^{\times}_{K_{s,y_{i'}}}$ holds simultaneously with the inclusion ${\left(\widetilde{\alpha}^{\ell_0}\right)}^{1/p_2} \in \oo^{\times}_{K_{w,i'}}$. Then the equality $p = p_2$ must hold in order for $K_{s, y_{i'}} \left( {\left(\widetilde{\alpha}^{\ell_0}\right)}^{1/p_2} \right)$ to be an intermediate field of the field extension $K_{w, i'} / K_{s, y_{i'}}$. However, then we have arrived at a contradiction because $p$ has maximal index in $P$, but the fact that $p = p_2$ implies that $p$ is succeeded by $p_1$ in $P$.
    \smallskip
\end{proof}

%%%%%%%%%%%%%%%
\bigskip

\section{Atomicity in Monoid Algebras of Exponentially Cyclic Puiseux Monoids} \label{sec: atomic exp cyc}
In this section, we first prove a key proposition governing the number of divisors of irreducible polynomials in $\zz[x]$ in certain subdomains of $\zz[M_q]$. Then we prove Theorem~\ref{thm: Z[M_q] is atomic}, our main result on the atomicity of the monoid algebra $\qq[M_q]$.

\begin{prop} \label{prop: (general) num irr. div is limited}
    Let $r \in \nn_{\geq 3}$ be an odd squarefree positive integer such that $\pi(r) = \{p_1, \ldots, p_m \}$ for some index $m \in \nn$ and odd primes $p_1,\ldots,p_m$. For brevity, set $\pi := \pi(r)$. Furthermore, let $p_{\text{max}}$ be the maximal prime in $\pi$. Let $f(x)$ be a nonconstant irreducible polynomial in $\zz[x] \setminus (\{ \pm x \} \cup \{ \pm \Phi_y(x): y \in \nn \})$ of degree $d \in \nn$. Then for all $\lambda \in \nn_0$, the polynomial $f(x)$ has at most 
    \[
    A(f(x),\pi) :=
    \begin{dcases}
        {r}^{
            \left\lceil
            \log_{\frac{p_{\text{max}}-1}{p_{\text{max}}}}\left(\log_{|\lc(f)|}{2}\right)\right\rceil
                + 1 + m \Lambda^{\star}(f(x), \pi)
            },&
        \lc(f) \notin \{ \pm 1 \},\\
            {r}^{\Lambda^{\star}(f(x),\pi)}, &
            \lc(f) \in \{ \pm 1 \},
    \end{dcases}
    \]
    irreducible divisors in the monoid algebra $\zz \left[\frac{1}{r^{\lambda}} \nn_0 \right]$. In particular, $A(f(x), \pi)$ is finite and well-defined by the assumption that $f(x)$ is not contained in the set $\{ \pm x \} \cup \{ \pm \Phi_y(x): y \in \nn \}$.
\end{prop}
\begin{proof}
    By isomorphism of algebras, this is equivalent to proving that $f(x^{r^{\lambda}})$ has at most $A(f(x),\pi)$ irreducible divisors in $\zz[x]$. We first take the case when $\lc(f) \notin \{ \pm 1 \}$. Set
    \[
     H := 
     \left\lceil 
     \log_{\frac{p_{\text{max}}-1}{p_{\text{max}}}}\left(\log_{|\lc(f)|}{2}\right)
     \right\rceil + 1
     ,
    \]
    $\theta_1 = \min(\lambda, H),$ and $\theta_2 := \lambda - \theta_1$.
    In light of Lemma~\ref{lemma: irr div appears in some splitting seq}, it is also sufficient to prove that for any splitting sequence, $F := (f_n(x))_{n \in \nn_0}$, of $f(x)$ with exponent sequence
    \[
    E_F := (r)_{n \in \ldb 1, \theta_1 \rdb} \circ (p_1)_{n \in \ldb 1, \theta_2 \rdb} \circ \cdots \circ (p_{m})_{n \in \ldb 1, \theta_2 \rdb},
    \]
    the polynomial $f_{\theta_1 + m \theta_2}(x)$ has degree at least $\frac{d r^{\lambda}}{A(f(x), \pi)}$. If $\lambda \leq H$, then $f(x^{r^{\lambda}})$ has at most $r^H$ irreducible divisors in $\zz[x]$ by successive applications of Lemma~\ref{lemma: <= r irr. divs}. In a similar fashion, if $\lambda$ is greater than $H$ but the inequality $\textsf{N}_{\text{spl}}(F) \leq H$ holds, then the previous statement also holds by Corollary~\ref{corll: SSS...}. In both cases, we are done, so henceforth assume that the quantities $\lambda$ and $\textsf{N}_{\text{spl}}(F)$ are both larger than $H$. Observing that
    \[
    {|\lc(f)|}^{\prod_{i = 1}^{\theta_1}{\frac{p_{\text{max}} - 1}{p_{\text{max}}}}} < 2
    \]
    by the definition of $H$, we can deduce that the leading coefficient of $f_{\theta_1}(x)$ is in the set $\{ \pm 1\}$ by Lemma~\ref{lemma: |lc| <= |lc|}.
    Observe that for every $i \in \ldb 0, m-1 \rdb$, the inclusion $f_{\theta_1+i \theta_2}(x) \in \mathcal{D}(f(x), \pi)$ holds.  Moreover, the subsequence 
    \[
    F_i := (f_n(x))_{n \in \ldb \theta_1+i \theta_2, \theta_1+(i+1) \theta_2 \rdb}
    \]
    is a $(p_{i+1})_{n \in \ldb 1, \theta_2 \rdb}$ subsplitting sequence of $F$. Therefore, Corollary~\ref{corll: SSS...}, in tandem with Lemma~\ref{lemma: Lambda(b(x),P) <= Lambda^*(f(x),P)}, guarantees that the total number of splits which occur in $F_i$ is bounded above by $\Lambda^{\star}(f(x), \pi)$. Therefore, the total number of splits in $F$ is at most $H + m \Lambda^{\star}(f(x), \pi)$.
    Since each split reduces the degree of the subsequent polynomial in $F$ by a factor of $r$ in the worst case (this is a direct consequence of Lemma~\ref{lemma: <= r irr. divs}), we have that
    \[
    \deg f_{\theta_1+ m \theta_2} \geq \frac{d r^{\lambda}}{
    {r}^{H + m \Lambda^{\star}(f(x), \pi)}
    }
    \]
    as required. A similar proof structure may be used for the case when $\lc(f) \in \{ \pm 1 \}$, although the quantities $H$ and $\theta_1$ are not needed for the proof.
\end{proof}

\begin{thm} \label{thm: Z[M_q] is atomic}
    The monoid algebra $\qq[M_q]$ is atomic for all $q \in ((0, 1) \cap \qq) \setminus \nn^{-1}_{> 1}$ with an odd denominator.  
\end{thm}
\begin{proof}
    We follow the proof of~\cite[Proposition~4.8]{BGLZ24}. First let $q \in ((0, 1) \cap \qq) \setminus \nn^{-1}_{> 1}$ such that $\textsf{d}(q)$ is odd and let $r$ be the squarefree part of $\textsf{d}(q)$. To show that $\qq[M_q]$ is atomic, Proposition~\ref{prop: Z[M] atomic = Q[M] atomic} guarantees that we need only prove the atomicity of $\zz[M_q]$. The claim is trivially true for any constant polynomial in $\zz[M_q]$. Accordingly, select any nonconstant polynomial $f(x) \in \zz[M_q]$. Since $M_q$ satisfies the almost ACCP, there exists some common divisor $d \in M_q$ of $\supp f$ such that $f(x)$ can be decomposed as $f(x) = x^d \cdot \frac{f(x)}{x^d}$, where the latter term is ACCP-supported. As $x^d$ is a monomial, it is atomic by virtue of the atomicity of $M_q$. Then the proof of the theorem would follow from proving the atomicity of the polynomial $\frac{f(x)}{x^d}$. To do this, we show that every ACCP-supported polynomial in $\zz[M_q]$ satisfies the ACCP.

    Set $R := \zz[M_q]$. The claim is trivially true for monomials so we can freely consider only the case of non-monomials. Accordingly, suppose by way of contradiction that there exists a polynomial $P(x) \in R$ with $|\supp P| > 1$ such that $P(x)$ is ACCP-supported but does not satisfy the ACCP in $R$. Then, setting $P_0(x) := P(x)$, there exists a sequence of $R$-ideals $\ppp := (P_n(x))_{n \in \nn_0}$ such that 
    \begin{equation} \label{eq: asc chain P}
        P_0(x) R \subsetneq P_1(x) R \subsetneq P_2(x) R \subsetneq  \cdots
    \end{equation}
    is a strictly ascending chain of principal ideals which does not stabilize.
    For each $n \in \nn_0$, let $Q_{n+1}(x) := \frac{P_{n}(x)}{P_{n+1}(x)}$ and set $Q := (Q_n(x))_{n \in \nn}$. Note that the strictness of~\eqref{eq: asc chain P} ensures that $\pm 1 \notin Q$. Now observe that only finitely many terms in $Q$ can have nonzero order as $P_0(x)$ is ACCP-supported. Moreover, since $\zz$ is Noetherian, only finitely many terms in $Q$ can have nonunit content. Then, by dropping finitely many initial terms in $\ppp$, we may pass to a subsequence of $\ppp$ (and hence to a subsequence of $Q$) such that the equalities $\ord Q_n = 0$ and $|\textsf{c}(Q_n)| = 1$ hold for all $n \in \nn$. 
    \smallskip

    \noindent \textsc{Claim.} 
    Let $\ell$ be an index in $\nn_0$ and $f(x)$ an irreducible polynomial in the monoid algebra $\zz[\frac{1}{r^{\ell}} \nn_0]$ such that $f(x^{r^{\ell}}) \notin \{ \pm x \} \cup \{ \pm \Phi_y(x): y \in \nn \}$. Then there exists $ \Gamma\left(f(x), \zz\left[\frac{1}{r^{\ell}} \nn_0\right]\right) \in \nn_0$ such that the polynomial $f(x)$ has the same factorization into irreducibles in the monoid algebra $\zz \left[\frac{1}{r^{\ell + \gamma}} \nn_0 \right]$ for every $\gamma \geq \Gamma\left(f(x), \zz\left[\frac{1}{r^{\ell}} \nn_0\right]\right)$.
    \smallskip

    \noindent \textsc{Proof of Claim.}
    If $f(x)$ is constant, the claim is trivial. So henceforth assume otherwise. As observed in Remark~\ref{rem: atoms of z[1/p^k N] = z[1/p^l N]}, it suffices to prove the claim for $\ell = 0$. Then $f(x)$ is an irreducible polynomial in $\zz[x]$ that is not contained in the set $\{ \pm x \} \cup \{ \pm \Phi_y(x): y \in \nn \}$: the claim follows directly from Proposition~\ref{prop: (general) num irr. div is limited}, in tandem with the fact that
    \[
    \zz[x] \subsetneq \zz\left[\frac{1}{r} \nn_0\right] \subsetneq \zz\left[\frac{1}{r^2} \nn_0\right] \subsetneq \zz\left[\frac{1}{r^3} \nn_0\right] \subsetneq \cdots
    \]
    is a tower of ring extensions of UFDs that preserves nonunits.
    \smallskip
    
    Let $\ell_P$ be the minimal index in $\nn_0$ such that the inclusion $P(x) \in \zz\left[\frac{1}{r^{\ell_P}} \nn_0\right]$ holds.  
    We can clearly express $P(x)$ as follows, for indices $m,k \in \nn_0$ (not both zero), an integer $\mu$, and nonconstant irreducible polynomials $p_1(x),\ldots,p_k(x) \in \zz[x] \setminus \{\pm \Phi_y(x): y \in \nn\}$:
    \begin{equation}\label{eq: P(x) = }
        P(x) = \mu \left(\prod_{i=1}^m{\Phi_{y_i}(x^{r^{- \ell_P}})}
        \right)
        \left(\prod_{j=1}^{k}{p_j(x^{r^{-\ell_P}})}\right).
    \end{equation}
    Then set
    \[
    \Gamma^{\star} := \max_{
        \substack{
        f(x) \mid_{\zz\left[\frac{1}{r^{\ell_P}} \nn_0\right]} P(x)\\
        f(x) \text{ irreducible}\\
        f(x^{r^{\ell_P}}) \notin \{ \pm x \} \cup \{\pm \Phi_y(x): y \in \nn\}\cup \zz \setminus \zz^{\times}
        }
    }{\Gamma\left(f(x), \zz\left[\frac{1}{r^{\ell_P}} \nn_0\right]\right)}.
    \]
    Since $Q$ has infinitely many terms, the equality $\inf \{\deg Q_{n}: n \in \nn \} = 0$ holds, and so we may fix an index $X^{\star}$ in $\nn$ such that the inequality $\deg Q_{\chi} < \frac{1}{r^{\Gamma^{\star}}}$ holds for all $\chi > X^{\star}$. For each $\chi \in \nn_{> X^{\star}}$, we let ${\ell}_{\chi}$ be the minimal index in $\nn$ such that the inclusions $Q_{\chi}(x), \frac{P(x)}{Q_{\chi}(x)}, P(x) \in \zz\left[ \frac{1}{r^{{\ell}_{\chi}}} \nn_0 \right]$ hold. Then clearly $l_{\chi} \geq \ell_P$ and so $\zz\left[ \frac{1}{r^{\ell_P}} \nn_0 \right] \subseteq \zz\left[ \frac{1}{r^{{\ell}_{\chi}}} \nn_0 \right]$ is a ring extension of UFDs preserving nonunits. Therefore, for any irreducible divisor $q(x) \in \zz\left[ \frac{1}{r^{{\ell}_{\chi}}} \nn_0 \right]$ of $Q_{\chi}(x)$, there exists an irreducible divisor $f_{q}(x) \in \zz\left[ \frac{1}{r^{\ell_P}} \nn_0 \right]$ of $P(x)$ such that the divisibility relation $q(x) \mid_{\zz\left[ \frac{1}{r^{{\ell}_{\chi}}} \nn_0 \right]} f_q(x)$ holds. Then, for any $\chi \in \nn_{> X^{\star}}$, the fact that the inequality $\deg Q_{\chi} < \frac{1}{r^{\Gamma^{\star}}}$ holds guarantees that $\deg q < \frac{1}{r^{\Gamma^{\star}}}$ and so $f_q(x^{r^{\ell_P}})$ must be in the set $\{ \pm x \} \cup \{ \pm \Phi_y(x): y \in \nn \} \cup \zz \setminus \zz^{\times}$. Noting that $\ord Q_{\chi} = 0$ and that $\textsf{c}(Q_{\chi}) \in \{\pm 1\}$, we may apply Lemma~\ref{lemma: cyclotomicFactors} to obtain the following divisibility relation, where $\iota(\chi) \in \nn$ is the number of irreducible divisors of $Q_{\chi}(x)$ in $\zz\left[ \frac{1}{r^{\ell_{\chi}}} \nn_0 \right]$ and the indices $y_{\chi,1}, \ldots, y_{\chi, \iota(\chi)}$ are (not necessarily distinct) natural numbers:
    \[
    Q_{\chi}(x) \mid_{\zz\left[ \frac{1}{r^{\ell_{\chi}}} \nn_0 \right]} \prod_{i=1}^{\iota(\chi)}{\Phi_{y_{\chi, i}}(x^{r^{- \ell_P}})}.
    \]
    Once again applying the cyclotomic identity given in Lemma~\ref{lemma: cyclotomicFactors}, in addition to a translation of factorizations between isomorphic monoid algebras (as in Remark~\ref{rem: atoms of z[1/p^k N] = z[1/p^l N]}), we conclude that there exists an index $h_{\chi} \in \nn$ and indices $w_{\chi,1},\ldots,w_{\chi, h_{\chi}} \in \nn$ such that for any $\chi \in \nn_{> X^{\star}}$, the polynomial $Q_{\chi}(x)$ factorizes as
    \[
    Q_{\chi}(x) = \prod_{i' = 1}^{h_{\chi}}{
    \Phi_{w_{\chi, i'}}(
    x^{r^{- {\ell}_{\chi}}}
    )}
    \]
    in $\zz\left[\frac{1}{r^{{\ell}_{\chi}}} \nn_0\right]$. Moreover, the fact that $\gcd(\textsf{n}(q),r) = 1$ implies that
    \begin{equation}\label{eq: v n(q)- adic ineq}
        \vp_{\textsf{n}(q)}\left(  \lcm\left(w_{\chi,1},\ldots,w_{\chi, h_{\chi}}\right)
        \right) =
        \vp_{\textsf{n}(q)}\left(\lcm\left(y_1,\ldots,y_m\right)
        \right).
    \end{equation}
    Observe that by Lemma~\ref{lemma: factor-x**n-1}, the following divisibility relation holds in $\zz[x]$:
    \begin{equation}\label{eq: x**u - 1}
        \prod_{i' = 1}^{h_{\chi}}{
        \Phi_{w_{\chi, i'}}(x)} \mid_{\zz[x]} x^{\lcm\left(w_{\chi,1},\ldots,w_{\chi, h_{\chi}}\right)} - 1.
    \end{equation}
    For each $\chi \in \nn_{> X^{\star}}$, set $V_{\chi} := \vp_{\textsf{n}(q)}\left(   \lcm\left(w_{\chi,1},\ldots,w_{\chi, h_{\chi}}\right)
    \right)$. The equality $\inf \{\deg Q_{n}: n \in \nn \} = 0$, in conjunction with~\eqref{eq: v n(q)- adic ineq}, ensures that we can find $\chi' \in \nn_{> X^{\star}}$ such that 
    $\deg Q_{\chi'} < q^{V_{\chi'}}$.
    Setting $c := \min(\supp Q_{\chi'} \setminus \{0\})$, the previous inequality guarantees that $c \in \left \langle q^v: v > V_{\chi'} \right \rangle \subseteq M_{q}$ and hence that $\vp_{\textsf{n}(q)}(c) > V_{\chi'}$. On the other hand, using~\eqref{eq: x**u - 1}, we may follow the proof of~\cite[Proposition~4.8]{BGLZ24} to see that there exist indices $t_1, t_2 \in \nn_0$ such that $c r^{t_1} \mid_{\zz}  r^{t_2} \lcm\left(w_{\chi',1},\ldots,w_{\chi', h_{\chi'}}\right)$. However, this implies that the $\textsf{n}(q)$-adic valuation of $c$ is at most $V_{\chi'}$, a contradiction. As such, we conclude that all ACCP-supported non-monomials in $\zz[M_q]$ satisfy the ACCP and so are atomic.
\end{proof}

\bigskip
%%%%%%%%%%%%%%
%%%%%%%%%%%%%%

\end{document}